\documentclass[10pt]{article}

\usepackage[T1]{fontenc}
\usepackage[utf8]{inputenc}
\usepackage{lmodern}
\usepackage{amsmath,amssymb,amsthm,mathtools}
\usepackage{booktabs,tabularx,array,placeins}
\usepackage{microtype}
\usepackage[margin=1.15in]{geometry}
\usepackage{xcolor}
\usepackage[hidelinks]{hyperref}
\hypersetup{pdftitle={Fast Computation of the Kippenhahn Hypersurface of a Joint Numerical Range},
 pdfauthor={Carlos Rito},
 pdfkeywords={joint numerical range, Kippenhahn variety, projective duality, modular interpolation}}

\newcommand{\C}{\mathbb{C}}
\newcommand{\R}{\mathbb{R}}
\newcommand{\Pj}{\mathbb{P}}

\newcommand{\tr}{\operatorname{tr}}
\newcommand{\conv}{\operatorname{conv}}
\newcommand{\reg}{\mathrm{reg}}

\newtheorem{proposition}{Proposition}[section]

\title{Fast Computation of the Kippenhahn Hypersurface\\ of a Joint Numerical Range}
\author{Carlos Rito}
\date{}

\begin{document}
\maketitle

\begin{abstract}
We present a fast modular interpolation method for computing the hypersurface
components of algebraic Kippenhahn varieties over number fields.
Starting from Hermitian matrices, we factor their determinantal
polynomial and use Katz's test to identify the factors whose duals are
hypersurfaces. For each such factor, we sample smooth points over finite
fields and compute their tangent hyperplanes. These give points on the
dual hypersurface, from which we interpolate its defining equation.
For number-field coefficients, we use primes at which the defining
polynomial of the field factors completely into distinct linear factors.
Interpolation across all its roots recovers the residues of the
coefficients in a fixed number-field basis. Chinese remaindering and
rational reconstruction then recover an equation over the original
number field from the computations at different primes.
An exact polynomial divisibility test and an irreducibility check certify
each reconstructed equation. Four determinantal quartics arising from
matrix examples in the literature yield certified dual surfaces of
degrees $12$, $18$, $24$ and $34$. We give implementations in Magma and
Singular and compare their running times with Gr\"obner basis elimination.

\end{abstract}

\section{Introduction}

The numerical range of a complex $d\times d$ matrix $A$ is
\[
 W(A)=\{x^*Ax:x\in\C^d,\ x^*x=1\}\subset\C.
\]
More generally, Hermitian matrices $A_1,\ldots,A_n$ determine the joint
numerical range
\[
 W(\mathbf A)=\{(x^*A_1x,\ldots,x^*A_nx):x\in\C^d,\ x^*x=1\}
 \subset\R^n.
\]
These sets have attracted sustained
interest in matrix and operator theory, numerical analysis, optimization,
convex and algebraic geometry, and quantum information.

In matrix analysis, the numerical range provides information about a
matrix that its eigenvalues alone do not capture. The theorem of
Crouzeix and Palencia \cite{CrouzeixPalencia2017} bounds the norm of an
analytic function of a matrix in terms of the values of that function
on the numerical range. In numerical linear algebra, Beckermann,
Goreinov and Tyrtyshnikov \cite{BeckermannGoreinovTyrtyshnikov2006} use
its geometry to sharpen convergence estimates for the generalized
minimal residual method (GMRES), an iterative method for solving
nonsymmetric linear systems. Related convexity questions also arise
in optimization and control: Polyak
\cite{Polyak1998} studies quadratic transformations and applies their
convexity properties to quadratic optimization and control problems.

In quantum information, joint numerical ranges describe the
simultaneously attainable expectation values of observables.
Xie et al.\ \cite{XieEtAl2020} observe joint
numerical ranges in a three-level photonic system and relate their
geometry to ground-state degeneracies.

The geometry of numerical ranges has itself generated a substantial
literature. Li and Poon \cite{LiPoon2000} and Gutkin, Jonckheere and
Karow \cite{GutkinJonckheereKarow2004} investigate convexity of joint
numerical ranges. Szyma\'nski, Weis and \.{Z}yczkowski
\cite{SzymanskiWeisZyczkowski2018} classify the joint numerical ranges
of three Hermitian matrices of size three, while Pikul et al.\
\cite{PikulEtAl2026} study the corresponding geometry for matrices of
size four. Chien and Nakazato \cite{ChienNakazato2002} use
Davis--Wielandt shells to describe the boundaries of $q$-numerical
ranges and derive a plotting algorithm. They also study flat portions
on the boundary of these shells for $3\times3$ matrices
\cite{ChienNakazato2009}.

An algebraic approach goes back to Kippenhahn,
who describes
the numerical range as the convex hull of a real algebraic curve obtained
by projective duality (see the English translation
\cite{ZachlinHochstenbach2008}).
Chien and Nakazato \cite{ChienNakazato2010} investigate generating
hypersurfaces for joint numerical ranges, and Plaumann, Sinn and Weis
\cite{PlaumannSinnWeis2021} establish a higher-dimensional form of
Kippenhahn's theorem, linking projective duality, convex geometry and
quantum states.

Numerical approximation of these sets can be obtained by solving
Hermitian eigenvalue problems in sampled directions. In the classical
case, algorithms for tracing the boundary include those of Johnson
\cite{Johnson1978} and Loisel and Maxwell \cite{LoiselMaxwell2018}.
Such methods give accurate geometric approximations, but do not
determine exact algebraic degrees, irreducible components or
singularities. These questions require equations for the associated
Kippenhahn varieties, whose exact computation remains costly.
Methods based on Gr\"obner basis elimination or resultants can be very
slow even for matrices of size four, as illustrated by the
Davis--Wielandt examples of Chien and Nakazato
\cite{ChienNakazato2025DW} and the comparisons in
Section~\ref{sec:timings}.

Our goal is to provide a fast method for computing these equations,
with exact certification, when the dual components are hypersurfaces.
We sample smooth points on the
determinantal components over finite fields and compute their tangent
hyperplanes, thereby obtaining points on the duals. Interpolation then
recovers an equation from these points. For coefficients in a number
field, we use primes at which its defining polynomial factors completely
and combine the equations obtained from all its roots before lifting
to characteristic zero. Exact divisibility and irreducibility tests
certify the resulting equations.
The algorithm computes equations for the algebraic dual hypersurfaces.
Extracting the relevant real sets and forming their convex hull, as
described in Section~\ref{sec:kippenhahn-theorem}, are separate tasks.

Section~\ref{sec:background} gives the definitions and the relevant
duality results. Section~\ref{sec:algorithm} describes the modular
algorithm and its exact certification. Section~\ref{sec:applications}
revisits four matrix examples from Chien and Nakazato
\cite{ChienNakazato2025DW} and Pikul et al.\ \cite{PikulEtAl2026},
recovers their dual hypersurface equations, and compares running times
for modular interpolation and Gr\"obner basis elimination in Magma and
Singular.

\bigskip
\noindent\textbf{Use of AI.}
The original algorithm was written in Magma by the author.
To provide an implementation in open-source software, the Magma code
was translated into Singular by ChatGPT~6.

\medskip
\noindent\textbf{Acknowledgments.}
The author thanks Gra\c ca Soares for introducing him to numerical ranges
through a talk on her work \cite{BebianoLemosSoares2023}, which sparked
his interest in the subject.
The author also thanks Hiroshi Nakazato for helpful email correspondence
about the computations in \cite{ChienNakazato2025DW}.
The author was financed by Portuguese funds through FCT
(Funda\c c\~ao para a Ci\^encia e a Tecnologia).  This work was carried out
within project UID/00013/2025 of the Centro de Matem\'atica da Universidade
do Minho (CMAT/UM).

\section{Numerical ranges and Kippenhahn varieties}
\label{sec:background}

This section recalls the definitions and projective-duality results that
underlie the computation of Kippenhahn varieties from joint numerical
ranges.

\subsection{Numerical and joint numerical ranges}

Let $\mathcal H_d$ be the real vector space of $d\times d$ Hermitian
matrices, and fix $A_1,\ldots,A_n\in\mathcal H_d$.
A \emph{density matrix} is a matrix $\rho\in\mathcal H_d$ that is
positive semidefinite and has trace one. We write $\rho\succeq0$ for
positive semidefiniteness, meaning that $v^*\rho v\geq0$ for every
$v\in\C^d$.

A \emph{pure state} is a density matrix of rank one. Such a matrix has
the form $\rho=xx^*$ for a unit vector $x\in\C^d$, and
\[
 \tr(\rho A_i)=\tr(xx^*A_i)=x^*A_ix.
\]
Thus the joint numerical range $W(\mathbf A)$ defined in the introduction
is obtained by allowing $\rho$ to range over the pure states.
Allowing all density matrices gives the set
\[
 \widehat W(\mathbf A)
 =\{(\tr(\rho A_1),\ldots,\tr(\rho A_n)):
       \rho\in\mathcal H_d,\ \rho\succeq0,\ \tr(\rho)=1\}.
\]
By the spectral theorem, every density matrix has a decomposition
\[
 \rho=\sum_{j=1}^d t_jx_jx_j^*,
\]
where $x_1,\ldots,x_d\in\C^d$ form an orthonormal basis of eigenvectors
of $\rho$, and $t_1,\ldots,t_d$ are the corresponding eigenvalues.
Since $\rho$ is positive semidefinite and has trace one, $t_j\geq0$
and $\sum_{j=1}^d t_j=1$. Hence
\[
 \tr(\rho A_i)=\sum_{j=1}^d t_j\,x_j^*A_ix_j,
 \qquad i=1,\ldots,n,
\]
so every point of $\widehat W(\mathbf A)$ is a convex combination of
points of $W(\mathbf A)$. Conversely, every convex combination of
rank-one matrices $xx^*$ with $\|x\|=1$ is a density matrix. Therefore
\[
 \widehat W(\mathbf A)=\conv W(\mathbf A).
\]
This convex set is also called the convex support or the joint algebraic
numerical range.

For two Hermitian matrices, the Toeplitz--Hausdorff theorem
(see \cite[Section~1]{PlaumannSinnWeis2021}) gives
\[
 W(A_1,A_2)=\widehat W(A_1,A_2).
\]
Under the identification $\R^2\cong\C$ given by
$(y_1,y_2)\mapsto y_1+\mathrm{i}y_2$, this is the usual numerical range
$W(A_1+\mathrm{i}A_2)$.
The joint numerical range of three Hermitian matrices is also convex
when $d\geq3$, by a theorem of Au-Yeung and Poon
\cite[Theorem~2]{AuYeungPoon1979}, see also
\cite{GutkinJonckheereKarow2004}.
Consequently, $W=\widehat W$ in all the $4\times4$ applications of
Section~\ref{sec:applications}.

\subsection{Kippenhahn's theorem and its higher-dimensional form}
\label{sec:kippenhahn-theorem}

The relation between numerical ranges and projective duality originates
in Kippenhahn's theorem (see \cite{ZachlinHochstenbach2008}). Let
$A_1,A_2\in\mathcal H_d$ and consider the homogeneous determinantal
polynomial
\[
   p_{\mathbf A}(u_0,u_1,u_2)
   =
   \det(u_0I+u_1A_1+u_2A_2).
\]
It defines a complex projective plane curve
\[
   X_{\mathbf A}=V(p_{\mathbf A})\subset\Pj^2_\C.
\]
Let $X_{\mathbf A}^{\vee}\subset(\Pj^2_\C)^{\vee}$ denote its projective
dual variety, obtained from the tangent lines to the smooth loci of
its reduced components. Linear components contribute points. We view
its real affine part, in the chart where the first dual coordinate is
one, as a subset of the real plane:
\[
   \mathcal K(\mathbf A)
   =
   \left\{
      (y_1,y_2)\in\R^2
      \;:\;
      (1:y_1:y_2)\in X_{\mathbf A}^{\vee}(\R)
   \right\}.
\]
This is usually called the \emph{Kippenhahn curve} or the
\emph{boundary-generating curve}. Kippenhahn's theorem states that
\[
   W(A_1,A_2)
   =
   \operatorname{conv}\mathcal K(\mathbf A).
\]
Thus, although the boundary of the numerical range need not itself be
algebraic, the entire numerical range is generated, by convexification,
from an algebraic curve obtained by projective duality.

This result has a higher-dimensional analogue. For
$\mathbf A=(A_1,\ldots,A_n)$, consider
\[
   p_{\mathbf A}(u_0,\ldots,u_n)
   =
   \det\left(u_0I+\sum_{j=1}^n u_jA_j\right)
\]
and the determinantal hypersurface
\[
   X_{\mathbf A}=V(p_{\mathbf A})\subset\Pj^n.
\]
All varieties are understood with their reduced structure: if
\[
   p_{\mathbf A}=\prod_{i=1}^r f_i^{m_i}
\]
is its factorization into distinct geometrically irreducible factors, we set
$X_i=V(f_i)$ and disregard the multiplicities $m_i$. The union
\[
   X_{\mathbf A}^{\vee}
   =
   X_1^{\vee}\cup\cdots\cup X_r^{\vee}
   \subset(\Pj^n)^{\vee}
\]
will be called the \emph{algebraic Kippenhahn variety} associated with
$\mathbf A$.

To relate this projective variety to the numerical range in $\R^n$, we
use the affine chart $y_0\ne0$ in the dual projective space and normalize
$y_0=1$. This choice reflects the coefficient $I$ of $u_0$ in
$u_0I+u_1A_1+\cdots+u_nA_n$: a density matrix $\rho$ determines a point
of the ambient dual projective space with coordinates
\[
   (\tr(\rho I):\tr(\rho A_1):\cdots:\tr(\rho A_n))
   =
   (1:\tr(\rho A_1):\cdots:\tr(\rho A_n)),
\]
since $\tr(\rho)=1$. The remaining affine coordinates are precisely
those used to define $\widehat W(\mathbf A)$.

In higher dimensions, not every real singular point of
$X_{\mathbf A}^{\vee}$ is relevant to the joint numerical range. One
considers
\[
   T_i
   =
   \left\{
      (y_1,\ldots,y_n)\in\R^n
      \;:\;
      (1:y_1:\cdots:y_n)\in
      (X_i^{\vee})_{\mathrm{reg}}
   \right\}
\]
and the Euclidean closure
\[
   \mathcal K_{\R}(\mathbf A)
   =
   \overline{T_1\cup\cdots\cup T_r}.
\]
Here $(X_i^{\vee})_{\mathrm{reg}}$ denotes the smooth locus of
$X_i^{\vee}$. Taking the closure restores the singular points that are
limits of smooth real points, while excluding singular real points that
arise only from complex branches.

The generalized Kippenhahn theorem of Plaumann, Sinn and Weis
\cite{PlaumannSinnWeis2021} states that
\[
   \widehat W(\mathbf A)
   =
   \operatorname{conv}\mathcal K_{\R}(\mathbf A).
\]
Using our convention for the joint numerical range, this is equivalently
written as
\[
   \operatorname{conv}W(\mathbf A)
   =
   \operatorname{conv}\mathcal K_{\R}(\mathbf A).
\]
In particular, whenever $W(\mathbf A)$ is convex,
\[
   W(\mathbf A)
   =
   \operatorname{conv}\mathcal K_{\R}(\mathbf A).
\]
It follows that projective duals of the irreducible components of the
determinantal hypersurface are fundamental algebraic objects in the study
of joint numerical ranges. Their appropriate real parts contain the
extreme points needed to generate the convex hull of $W(\mathbf A)$.

\section{A modular algorithm for Kippenhahn varieties}
\label{sec:algorithm}

The input is a tuple $\mathbf A=(A_1,\ldots,A_n)$ of $s\times s$
Hermitian matrices over a number field $K$, with $n\geq2$. Our aim is
to compute equations for the hypersurface components of the algebraic
Kippenhahn variety. We first form and factor
\[
 p_{\mathbf A}=\det\left(u_0I+\sum_{j=1}^n u_jA_j\right)
              =\prod_{i=1}^r F_i^{m_i}
\]
over $K$, where the $F_i$ are distinct homogeneous irreducible factors.
The multiplicities can be ignored, since
$X_{\mathbf A}^{\vee}=\bigcup_{i=1}^r V(F_i)^{\vee}$.
We process each factor $F=F_i$ separately. Katz's test first determines
whether its dual is a hypersurface. If it is not, we record the marker $0$
for that factor.

For a factor that passes this test, we seek an irreducible homogeneous
polynomial $G(y_0,\ldots,y_n)$ defining $V(F)^{\vee}$. At a smooth point $x$ of
$V(F)$, the tangent hyperplane has coefficients $F_{u_i}(x)$, where
$F_{u_i}=\partial F/\partial u_i$. Thus
$(F_{u_0}(x):\cdots:F_{u_n}(x))$ is a point on the dual. Sampling smooth
points of $V(F)$ therefore supplies points at which $G$ must vanish.
For a fixed trial degree, these vanishing conditions are linear
conditions on the unknown coefficients of $G$, so a candidate equation
can be found by linear algebra.

To keep the arithmetic inexpensive, we perform this sampling and
interpolation over finite fields. Write
\[
 K=\mathbb Q(\alpha),
\]
where $\alpha$ has monic irreducible minimal polynomial
$f\in\mathbb Q[t]$. We choose primes $p$ for which $f$ reduces
to a product of distinct linear factors over $\mathbb F_p$.
Each root of this reduction gives a way to reduce the coefficients of
$F$ to $\mathbb F_p$, and we carry out the computation for every root.
For each resulting polynomial, we sample smooth points, compute their
tangent hyperplanes, and search for a degree at which these points
determine a unique homogeneous equation up to scalar. The degree is
estimated during the computation.

For each retained prime, interpolation across the roots of $f$ recovers
the residues of the rational coordinates of the unknown coefficients
in a power basis of $K$. We repeat the computation at further primes,
combine compatible modular data using the Chinese remainder theorem,
and apply rational reconstruction to obtain a candidate
$H\in K[y_0,\ldots,y_n]$ for the unknown equation $G$.
When this candidate remains unchanged after
adding another prime, we test exactly over $K$ that
$F\mid H(F_{u_0},\ldots,F_{u_n})$ and that $H$ is irreducible.
If both tests succeed, Proposition~\ref{prop:certificate} below shows that
$H$ is a nonzero scalar multiple of $G$. With both leading coefficients
normalized to one, this means $H=G$.
Repeating the procedure for all factors gives a list of
certified equations.

Irreducibility here is over $K$. If $F$ splits over an extension, its
returned equation defines the union of the duals of all geometric
components of $V(F)$. Sampling may miss some of these components, but
the exact divisibility test checks vanishing on every one of their
duals.

Throughout, we fix a degree-compatible monomial order and normalize
$F$, $G$ and each candidate equation to have leading coefficient one.
Before reduction, we exclude primes dividing denominators of
coefficients of $f$ or of the rational power-basis coordinates of the
coefficients of this normalized $F$. Infinitely many completely split
primes remain, by the Chebotarev density theorem applied to the splitting
field of $f$ over $\mathbb Q$. The following subsections give the details
of each stage.

\subsection{Katz's test}
\label{sec:katz-test}

Before interpolating, we determine whether $V(F)^{\vee}$ is a
hypersurface. Put $X=V(F)$, and let $\operatorname{Hess}(F)$ denote
the matrix of second partial derivatives of $F$. For $\deg F\geq2$,
let $h$ be the rank of the Hessian at a general point of $X$.
The hypersurface form of Katz's dimension formula
\cite[Chapter~1, Sections~5.A--B]{GelfandKapranovZelevinsky1994}
gives
\[
 \dim X^{\vee}=h-2.
\]
Since $\operatorname{Hess}(F)$ has size $(n+1)\times(n+1)$, the dual
$X^{\vee}$ is a hypersurface if and only if the Hessian determinant
does not vanish identically on $X$. By the irreducibility of $F$, this
is equivalent to
\begin{equation}\label{eq:katz-test}
 F\nmid\det\operatorname{Hess}(F).
\end{equation}

The function \texttt{KatzTest} computes the Hessian determinant and tests
its divisibility by $F$ over $K$. A nonzero remainder gives a positive
answer. A zero remainder means that the dual has smaller dimension.

The main routine records $0$ for every negative answer and
applies interpolation and exact certification to each factor that
passes the test.

\bigskip
\noindent\emph{Remark.} The algorithm can be adapted to compute defining
equations for dual varieties of higher codimension. We restrict ourselves
to the hypersurface case to keep the paper short and simple.

\subsection{The Gauss map and interpolation}

Let $X=V(F)$. The Gauss map on the smooth locus is
\[
 \gamma_F:X_{\reg}\longrightarrow(\Pj^n)^{\vee},\qquad
 x\longmapsto(F_{u_0}(x):\cdots:F_{u_n}(x)).
\]
Its image closure is $X^{\vee}$. Assume this is a hypersurface, with
irreducible equation $G\in K[y_0,\ldots,y_n]$ of degree $D$, normalized
to have leading coefficient one.
For $\delta=\deg F$, it is known that
\begin{equation}\label{eq:degree-bound}
 D\leq\delta(\delta-1)^{n-1}.
\end{equation}
See \cite[Lemma~1]{Ballico2003}.

Write $m=[K:\mathbb Q]$. For each admissible prime $p$, the reduction
$\bar f$ of $f$ modulo $p$ factors as
\[
 \bar f(t)=\prod_{j=1}^m(t-r_j),\qquad r_j\in\mathbb F_p,
\]
with distinct roots. Express each coefficient of $F$ in the power basis
$1,\alpha,\ldots,\alpha^{m-1}$, reduce its rational coordinates modulo
$p$, and substitute $\alpha=r_j$. This gives a polynomial over
$\mathbb F_p$ for each root.

We apply the following procedure to each resulting polynomial, denoted
by $F_p$. We choose a variable $u_\ell$ in which $F_p$ has the smallest positive
degree and assign random elements of the field $\mathbb F_p$ to the remaining
variables. Whenever the resulting univariate polynomial is nonconstant
and has roots in $\mathbb F_p$, we choose one root at random. Roots from
the same specialization can yield linearly dependent conditions in the
subsequent interpolation. Keeping only one root aims to avoid enlarging
the evaluation matrix with such redundant conditions.
This gives a point $x\in V(F_p)$.
If $\nabla F_p(x)\ne0$, its Gauss image lies on
the dual. Repeating this process and retaining distinct images gives
a set $S_p$ of sampled points on the dual.

To avoid an indefinite search, we abandon a prime after
$\max\{1000,100\deg_{u_\ell}F_p\}$ consecutive specializations fail
to produce a new Gauss image. The counter is reset whenever a new
image is found. If the cutoff is reached for any root of $\bar f$,
we discard the whole prime and continue with the next admissible prime.

We now seek a homogeneous equation through these points. For a trial
degree $e$, the Magma function \texttt{LinearSystem} finds all such
forms by computing the nullspace of the evaluation matrix $E_{p,e}$.
This matrix has one row for each point of $S_p$ and one column for each
degree-$e$ monomial. Its entries are the values of these monomials evaluated at
coordinates of the points.
The kernel
\[
 L_{p,e}=\ker E_{p,e}
\]
consists of the coefficient vectors of degree-$e$ forms vanishing
at every sampled point.
When it is one-dimensional, its generator gives a unique equation of
that degree through the sample, up to a scalar. We normalize it to have
leading coefficient one with respect to the fixed monomial order.

At primes where reduction preserves the dual hypersurface and its
degree, this recovers the corresponding reduction of $G$ when $e=D$.
It requires $\binom{n+D}{n}-1$ independent point conditions.
The next subsection explains how we choose the trial degrees and
the number of points to sample.

\subsection{Candidate degree selection}
\label{sec:degree-search}

We seek the degree for which \texttt{LinearSystem} returns exactly one
equation, up to a nonzero scalar. We begin with a trial degree $d_0$.
With sufficiently many general points, the outcome tells us how to
proceed: no equation means that the degree is too small, while several
independent equations mean that it is too large.

The initial upper bound of the search interval is
$\delta(\delta-1)^{n-1}$. After each trial, we search above
the tested degree if there is no equation, and below it if there is
more than one independent equation. The next trial is the midpoint
of the remaining interval, rounded up to an integer.
As soon as there is exactly one equation, we retain it together with
its degree.

For a trial degree $e$, we sample at least
\[
 \left\lceil\frac{11}{10}\binom{n+e}{n}\right\rceil
\]
distinct Gauss images, keeping all points already collected. This is a
$10\%$ surplus over the number of degree-$e$ monomials. The surplus is an
empirical choice intended to reduce the risk of an incorrect candidate
degree. It does not guarantee independent interpolation conditions or
the correct degree, and the user can easily increase it in the code if
needed.

The degree found is used as the first trial at the next root.
Computations for each next prime typically require only one interpolation.

\subsection{Reconstructing a candidate equation}
\label{sec:reconstruction}

Recall that our coefficients lie in the number field $K=\mathbb Q(\alpha)$,
where $\alpha$ has minimal polynomial $f$ of degree $m=[K:\mathbb Q]$.
We now recover an equation over $K$ from the equations computed over
finite fields. For each admissible prime $p$, the preceding computation
gives equations $H_{p,1},\ldots,H_{p,m}$ over $\mathbb F_p$, one for
each root $r_j$ of $\bar f$, each normalized to have leading coefficient
one. If their leading monomials differ, we discard the whole prime.
Across primes, we retain the largest leading monomial encountered.
A prime giving a smaller leading monomial is discarded. A larger one
causes us to discard all accumulated modular data and the previous
candidate, and to start a new reconstruction with the current prime.
Equations with the same leading monomial are combined.
Since they are homogeneous, they also have the same degree.

If a coefficient of the normalized desired polynomial $G$ has a rational power-basis
coordinate with denominator divisible by $p$, then at least one of the
modular equations has a leading monomial smaller than that of $G$.
This is handled by the preceding rule.

Fix a monomial and write its coefficient in $G$ as
\[
 c=a_0+a_1\alpha+\cdots+a_{m-1}\alpha^{m-1},
 \qquad a_i\in\mathbb Q.
\]
For a fixed prime $p$, let $c_j$ be the coefficient of that monomial
in $H_{p,j}$. When the common leading monomial is that of $G$, these are
the reductions of $c$ obtained by substituting $\alpha=r_j$.
Writing $\bar a_i$ for the reduction of $a_i$ modulo $p$, we obtain
the linear system over $\mathbb F_p$
\[
 \begin{pmatrix}
  1 & r_1 & \cdots & r_1^{m-1}\\
  1 & r_2 & \cdots & r_2^{m-1}\\
  \vdots & \vdots & & \vdots\\
  1 & r_m & \cdots & r_m^{m-1}
 \end{pmatrix}
 \begin{pmatrix}
  \bar a_0\\
  \bar a_1\\
  \vdots\\
  \bar a_{m-1}
 \end{pmatrix}
 =
 \begin{pmatrix}
  c_1\\
  c_2\\
  \vdots\\
  c_m
 \end{pmatrix}.
\]
This matrix is invertible because the roots $r_j$ are
distinct. Solving the system recovers the residues $\bar a_i$.

Repeating this at different primes gives residues of each rational
number $a_i$. For each $i$, the \emph{Chinese remainder theorem} (CRT)
combines these residues into one residue modulo $M$, the product of
the primes used.
We then use Magma's \texttt{RationalReconstruction} function or Singular's
\texttt{farey} function, which seeks a fraction with a small numerator and
denominator having the prescribed residue modulo $M$. Once $M$ is sufficiently large
relative to the numerators and denominators of the $a_i$, this recovers
these rational numbers and hence the coefficient $c\in K$.

Applying this procedure to every monomial gives a candidate equation
$H\in K[y_0,\ldots,y_n]$. If reconstruction fails or the candidate
changes, we continue with more primes. The implementation imposes no
fixed limit on the number of primes used for reconstruction.
Once the candidate remains unchanged after adding another prime, we
apply the exact certification described in the next subsection.

\subsection{Exact certification}
\label{sec:certification}

To certify that the reconstructed polynomial defines the dual hypersurface,
we use an exact polynomial divisibility test and an irreducibility test.

\begin{proposition}\label{prop:certificate}
Let $F\in K[u_0,\ldots,u_n]$ be homogeneous and irreducible over a field of
characteristic zero, and suppose that $V(F)^{\vee}$ is a hypersurface
with irreducible equation $G$. For a nonzero homogeneous polynomial
$H\in K[y_0,\ldots,y_n]$ of positive degree,
\[
 F\mid H(F_{u_0},\ldots,F_{u_n})
 \quad\Longleftrightarrow\quad G\mid H.
\]
In particular, this divisibility together with irreducibility of $H$
over $K$ certifies $H=cG$ for some $c\in K^*$.
\end{proposition}

\begin{proof}
The smooth locus of $V(F)$ is dense, and its Gauss image is dense in
$V(G)$. The composition $H(\nabla F)$ vanishes on this smooth locus
exactly when $H$ vanishes on $V(G)$. The respective vanishing ideals
are the prime ideals $(F)$ and $(G)$ over $K$, giving the equivalence.
If $H$ is irreducible and $G\mid H$, the two polynomials differ by a
nonzero scalar.
\end{proof}

Evaluating $H(F_{u_0},\ldots,F_{u_n})$ directly in Magma can be
expensive. We therefore use multivariate Horner's rule to make this
computation faster, followed by \texttt{IsDivisibleBy} to test
divisibility by $F$.

When $F$ and $H$ have rational coefficients, we scale them to primitive
polynomials with integer coefficients, meaning that the gcd of the
coefficients of each polynomial is $1$. Gauss's lemma then guarantees
that an integer
polynomial divisible by the primitive polynomial $F$ over $\mathbb Q$
has a quotient with integer coefficients. We can therefore perform
both the evaluation and the divisibility test using fast integer arithmetic.
For nonrational coefficients, we work over $K$.

If a reduction of $H$ modulo an admissible prime $p$ has the same degree
and is irreducible over $\mathbb F_p$, then $H$ is irreducible over $K$.
Since testing over finite fields is much faster, we first try at most three
distinct suitable primes. If none gives an irreducible reduction, we
test irreducibility directly over $K$.

\subsection{The algorithm}
\label{sec:algorithm-procedure}

The input is a tuple $\mathbf A$ of Hermitian matrices
of the same size over $K$, a lower bound $p_0$ for the prime search
and an initial guess $d_0$ for the dual degree. We begin with the
first admissible prime $p\geq p_0$. The procedure is:
\begin{enumerate}
 \item Form the determinantal polynomial $p_{\mathbf A}$ and factor it
       over $K$.
 \item Apply Katz's test from Section~\ref{sec:katz-test} to each distinct
       factor.
 \item For each factor that passes the test, compute many random smooth
       points over $\mathbb F_p$ at successive completely split primes $p$
       and use their gradients to obtain points on the dual hypersurface.
 \item Compute the candidate degree as in Section~\ref{sec:degree-search}.
 \item Apply the leading-monomial checks of Section~\ref{sec:reconstruction}.
       For each retained prime, interpolate the coefficients using all
       roots of the defining polynomial of $K$ modulo $p$ and reconstruct
       a candidate equation over $K$ from the accumulated data.
 \item When a candidate equation remains unchanged after adding one more
       prime, apply the divisibility and irreducibility tests of
       Proposition~\ref{prop:certificate}.
 \item Return a sequence with one entry for each distinct factor:
       its certified dual equation or the marker $0$.
\end{enumerate}
The main function is called as
\begin{verbatim}
DualFromMatrices(A,p_0,d_0);
\end{verbatim}
One can also use
\begin{verbatim}
DualHypersurface(F,p_0,d_0);
\end{verbatim}
which computes the dual equation of the hypersurface $\{F=0\}$.

Although we do not claim successful termination for every run, the method
succeeded in all our experiments, suggesting a high probability of success.

Reduction modulo a prime may change the degree of the dual hypersurface.
This can lead to an incorrect candidate degree or equation, as can
insufficiently general samples. In particular, the sampled
points may fail to represent all geometric components of $V(F)$ when
$F$ is irreducible over $K$ but not geometrically irreducible.
If final certification fails, the computation stops, and the user must
restart it from scratch with an initial prime larger than every prime
used for reconstruction in the failed run, discarding the accumulated
modular data. The error message reports the last such prime.

\section{Four computational examples}
\label{sec:applications}

We use four determinantal quartics from the matrix examples in
\cite{ChienNakazato2025DW,PikulEtAl2026}, whose dual surfaces have degrees
$12$, $18$, $24$ and $34$. They illustrate finite-field interpolation,
reconstruction and exact certification at different computational scales.

The usual method for computing these dual equations is Gr\"obner basis
elimination. For a homogeneous irreducible polynomial $F$ of degree
$\delta\geq2$ whose dual is a hypersurface, we form the ideal
\[
 I=\bigl(F(u),\,y_0-F_{u_0}(u),\ldots,y_n-F_{u_n}(u)\bigr)
 \subset K[u_0,\ldots,u_n,y_0,\ldots,y_n].
\]
Eliminating $u_0,\ldots,u_n$ gives
\[
 I\cap K[y_0,\ldots,y_n]=(G),
\]
where $G$ defines the dual hypersurface.

We compute this elimination ideal in Magma using
\texttt{EliminationIdeal(I,n+1)} with default algorithm selection.
We work over $\mathbb Q$ whenever $F$ has rational coefficients, and
over $K$ otherwise. Magma's default Gr\"obner basis computations over
$\mathbb Q$ and number fields use Faug\`ere's $F_4$ algorithm with
modular reconstruction. For positive-dimensional ideals, the default
strategy usually computes a basis in a graded order and then converts
it to an elimination order by a Gr\"obner walk.
For Singular, we use \texttt{eliminate} on the same ideal.
Apart from the simplest degree-$12$ case, only the Singular computation
of the degree-$18$ case finished, but it was much slower than modular
interpolation.
The remaining elimination computations did not finish within $24$ hours.

Chien and Nakazato \cite{ChienNakazato2025DW} obtain the degree-$34$
dual equation using resultants and interpolation. Hiroshi Nakazato
estimated, in personal communication, that the full computation would
take around $12$ hours on a personal computer, including preparation.
He also provided their Mathematica files, which we used to reproduce
their computation on our computer in $163$ minutes.
To make the comparison fairer, we did not exploit the symmetry of the
polynomial in this reproduction. The reproduction still uses
example-specific information, notably the known degree $34$ and the
choice of factors retained from the resultants. For this reason, we do
not include this comparison in the table in Section~\ref{sec:timings}.
Starting directly from the tuple of Hermitian matrices, with no
example-specific preparation, our method computes and certifies the
dual equation in less than 6 minutes on the setup described in
Section~\ref{sec:timings}.

\subsection{The four inputs}
\label{sec:DW-applications}

For $T\in M_4(\C)$, let $T^*$ denote its conjugate transpose. We use the
three Hermitian matrices
\[
 \operatorname{Re}T=\frac{T+T^*}{2},\qquad
 \operatorname{Im}T=\frac{T-T^*}{2\mathrm i},\qquad T^*T.
\]
The Davis--Wielandt shell is their joint numerical range
\[
 \operatorname{DW}(T)
 =W\bigl(\operatorname{Re}T,\operatorname{Im}T,T^*T\bigr)
 \subset\R^3.
\]
The associated determinantal quartic is
\[
 F_T(u_0,u_1,u_2,u_3)
 =\det\bigl(u_0I+u_1\operatorname{Re}T
                  +u_2\operatorname{Im}T+u_3T^*T\bigr).
\]
We consider
\[
 T_{\rm nil}=
 \begin{pmatrix}
 0&1&0&1\\0&0&1&0\\0&0&0&1\\0&0&0&0
 \end{pmatrix},\qquad
 T_{\rm CN}=
 \begin{pmatrix}
 3&1&1&0\\0&2&1&1\\0&0&-2&1\\0&0&0&-3
 \end{pmatrix}.
\]
Their dual surfaces have degrees $12$ and $34$, respectively.
The first matrix is $T_3$ from
\cite[Proposition~4.4]{ChienNakazato2025DW}. The second is the corrected
matrix from Proposition~3.1 of the same paper.
Hiroshi Nakazato has confirmed in personal communication that the bottom-right
matrix entry should be $-3$ rather than $3$, and that the coefficient
of $u_0^3u_3$ in the displayed polynomial $4F_T$ should be $124$
rather than $123$.

For an intermediate example, we use the exact Hermitian triple
$E_{12}=(B_1,B_2,B_3)$ from \cite[Section~4]{PikulEtAl2026}:
\[
 \begin{gathered}
 B_1=\begin{pmatrix}
 \sqrt3&0&0&-1\\
 0&-\sqrt3/2&0&1\\
 0&0&-\sqrt3/2&-\mathrm i\\
 -1&1&\mathrm i&0
 \end{pmatrix},\qquad
 B_2=\begin{pmatrix}
 0&0&0&\mathrm i\\
 0&3/2&0&-\mathrm i\\
 0&0&-3/2&1\\
 -\mathrm i&\mathrm i&1&0
 \end{pmatrix},\\[4pt]
 B_3=\operatorname{diag}(0,0,0,-2).
 \end{gathered}
\]
Its irreducible determinantal quartic has a dual of degree $18$ with
$691$ nonzero terms. The matrices have entries in
$\mathbb Q(\mathrm i,\sqrt3)$, and the determinantal polynomial has
nonrational coefficients. This example therefore tests interpolation
at the different roots of the defining polynomial before lifting.

The fourth input is the Hermitian triple denoted by $E_1$ in
\cite{PikulEtAl2026}:
\[
 A_1=\begin{pmatrix}
 0&\mathrm i&\mathrm i&-2\\
 -\mathrm i&0&\mathrm i&2\\
 -\mathrm i&-\mathrm i&0&2\mathrm i\\
 -2&2&-2\mathrm i&0
 \end{pmatrix},\qquad
 A_2=\begin{pmatrix}
 -2&1&1&-2\mathrm i\\
 1&1&1&2\mathrm i\\
 1&1&1&2\\
 2\mathrm i&-2\mathrm i&2&0
 \end{pmatrix},
 \qquad A_3=\operatorname{diag}(0,0,0,-4).
\]
The irreducible quartic
$\det(u_0I+u_1A_1+u_2A_2+u_3A_3)$ has a dual of degree $24$
with $2\,919$ nonzero terms. The three inputs other than
$E_{12}$ have entries in $\mathbb Q(\mathrm i)$. The example scripts
represent all four over the common field
\[
 K=\mathbb Q(\alpha),\qquad \alpha^4-\alpha^2+1=0,
 \qquad \alpha=e^{\pi\mathrm i/6},
\]
so that $\mathrm i=\alpha^3$ and $\sqrt3=2\alpha-\alpha^3$.
For $E_{12}$, interpolation uses all four distinct roots of this
defining polynomial at each selected prime. Rational input polynomials
are handled over $\mathbb Q$.

\subsection{Certified computations}
\label{sec:timings}

We use the four inputs in Section~\ref{sec:DW-applications}, with
lower bound $p_0=2003$ for the prime search, initial degree guess $d_0=18$,
random seed $1$ and a $10\%$ sampling surplus. For $E_{12}$, the first
completely split prime at or above $2003$ is $2017$. For the other three
examples, the computation takes place over $\mathbb Q$, and the first
prime used is $2003$.

We compare runs with one thread and ten threads. The parallelization comes from the native linear algebra routines
in Magma and Singular. Within each run, primes and their roots are processed
sequentially. We run the benchmarks one after another on the same computer.

\begin{table}[ht]
\centering
\small
\setlength{\tabcolsep}{3pt}
\begin{tabular}{c r r r r r r r c c}
\hline
 & & & & \multicolumn{6}{c}{Time (s)}\\
\cline{5-10}
 & & & & \multicolumn{2}{c}{Magma} & \multicolumn{2}{c}{Singular} & Magma & Singular\\
\cline{5-6}\cline{7-8}
input & degree & terms & primes & 1 thread & 10 threads & 1 thread & 10 threads & \shortstack{Gr\"obner\\elimination} & \shortstack{Gr\"obner\\elimination}\\
\hline
$T_{\rm nil}$ & 12 & 135 & 5 & 0.62 & 0.61 & 0.85 & 0.76 & 0.90 & 0.02\\
$E_{12}$ & 18 & 691 & 9 & 10.35 & 9.32 & 13.64 & 10.66 & $\dagger$ & 526.94\\
$E_1$ & 24 & 2919 & 18 & 30.65 & 24.48 & 46.57 & 27.2 & $\dagger$ & $\dagger$\\
$T_{\rm CN}$ & 34 & 4047 & 31 & 542.57 & 347.94 & 911.83 & 570.48 & $\dagger$ & $\dagger$\\
\hline
\end{tabular}
\caption{The symbol $\dagger$ means that elimination did not finish
within $24$ hours.}
\end{table}

\FloatBarrier
\noindent Times are elapsed seconds. The modular interpolation timings include
determinant construction, factorization, Katz's test, modular reconstruction
and exact certification. The Gr\"obner elimination attempts used
ten threads.
Exact certification accounted for approximately
$4\%$ to $30\%$ of the total modular interpolation time.

For $T_{\rm nil}$, $E_{12}$, $E_1$ and $T_{\rm CN}$, respectively,
peak resident memory usage (in MiB, with one/ten threads) was
$83/84$, $88/92$, $517/533$ and $1668/1697$ in Magma and
$40/41$, $56/55$, $257/253$ and $836/844$ in Singular.
Peak memory is measured for the whole process.

The computations were performed on an Intel Core i7-14700 with $15.2$ GiB RAM,
running Ubuntu 24.04.5 LTS, Magma V2.30-1 and Singular 4.3.2.

The ancillary files are provided in both Magma and Singular versions.
Each version consists of a file containing the main functions, a script
that loads this file and runs the four examples, a script for Gr\"obner
basis elimination, and a file containing the four computed dual polynomials.
The Singular code may not work properly on systems other than Linux.

\FloatBarrier
\begingroup
\small
\bibliographystyle{alpha}
\bibliography{References}

\newcommand{\etalchar}[1]{$^{#1}$}
\begin{thebibliography}{XZC{\etalchar{+}}20}

\bibitem[AYP79]{AuYeungPoon1979}
Yik-Hoi Au-Yeung and Yiu-Tung Poon.
\newblock A remark on the convexity and positive definiteness concerning
  {Hermitian} matrices.
\newblock {\em Southeast Asian Bulletin of Mathematics}, 3(2):85--92, 1979.

\bibitem[Bal03]{Ballico2003}
Edoardo Ballico.
\newblock An effective {Bertini} theorem over finite fields.
\newblock {\em Advances in Geometry}, 3(4):361--363, 2003.

\bibitem[BGT06]{BeckermannGoreinovTyrtyshnikov2006}
Bernhard Beckermann, Sergei~A. Goreinov, and Eugene~E. Tyrtyshnikov.
\newblock Some remarks on the {Elman} estimate for {GMRES}.
\newblock {\em SIAM Journal on Matrix Analysis and Applications},
  27(3):772--778, 2006.

\bibitem[BLS23]{BebianoLemosSoares2023}
N.~Bebiano, R.~Lemos, and G.~Soares.
\newblock On the numerical range of {Kac--Sylvester} matrices.
\newblock {\em Electronic Journal of Linear Algebra}, 39:242--259, 2023.

\bibitem[CN02]{ChienNakazato2002}
Mao-Ting Chien and Hiroshi Nakazato.
\newblock {Davis--Wielandt} shell and $q$-numerical range.
\newblock {\em Linear Algebra and its Applications}, 340(1--3):15--31, 2002.

\bibitem[CN09]{ChienNakazato2009}
Mao-Ting Chien and Hiroshi Nakazato.
\newblock Flat portions on the boundary of the {Davis--Wielandt} shell of
  3-by-3 matrices.
\newblock {\em Linear Algebra and its Applications}, 430(1):204--214, 2009.

\bibitem[CN10]{ChienNakazato2010}
Mao-Ting Chien and Hiroshi Nakazato.
\newblock Joint numerical range and its generating hypersurface.
\newblock {\em Linear Algebra and its Applications}, 432(1):173--179, 2010.

\bibitem[CN25]{ChienNakazato2025DW}
Mao-Ting Chien and Hiroshi Nakazato.
\newblock {Davis--Wielandt} shells of 4 by 4 matrices.
\newblock {\em Linear Algebra and its Applications}, 723:182--200, 2025.

\bibitem[CP17]{CrouzeixPalencia2017}
Michel Crouzeix and C{\'e}sar Palencia.
\newblock The numerical range is a $(1+\sqrt{2})$-spectral set.
\newblock {\em SIAM Journal on Matrix Analysis and Applications},
  38(2):649--655, 2017.

\bibitem[GJK04]{GutkinJonckheereKarow2004}
Eugene Gutkin, Edmond~A. Jonckheere, and Michael Karow.
\newblock Convexity of the joint numerical range: topological and differential
  geometric viewpoints.
\newblock {\em Linear Algebra and its Applications}, 376:143--171, 2004.

\bibitem[GKZ94]{GelfandKapranovZelevinsky1994}
Israel~M. Gelfand, Mikhail~M. Kapranov, and Andrei~V. Zelevinsky.
\newblock {\em Discriminants, Resultants, and Multidimensional Determinants}.
\newblock Birkh{\"a}user, Boston, 1994.

\bibitem[Joh78]{Johnson1978}
Charles~R. Johnson.
\newblock Numerical determination of the field of values of a general complex
  matrix.
\newblock {\em SIAM Journal on Numerical Analysis}, 15:595--602, 1978.

\bibitem[LM18]{LoiselMaxwell2018}
S{\'e}bastien Loisel and Peter Maxwell.
\newblock Path-following method to determine the field of values of a matrix
  with high accuracy.
\newblock {\em SIAM Journal on Matrix Analysis and Applications},
  39(4):1726--1749, 2018.

\bibitem[LP00]{LiPoon2000}
Chi-Kwong Li and Yiu-Tung Poon.
\newblock Convexity of the joint numerical range.
\newblock {\em SIAM Journal on Matrix Analysis and Applications},
  21(2):668--678, 2000.

\bibitem[Pol98]{Polyak1998}
Boris~T. Polyak.
\newblock Convexity of quadratic transformations and its use in control and
  optimization.
\newblock {\em Journal of Optimization Theory and Applications},
  99(3):553--583, 1998.

\bibitem[PSS{\etalchar{+}}26]{PikulEtAl2026}
Piotr Pikul, Ilya Spitkovsky, Konrad Szyma{\'n}ski, Stephan Weis, and Karol
  {\.Z}yczkowski.
\newblock Joint numerical ranges of three {Hermitian} {{\(4 \times 4\)}}
  matrices.
\newblock {\em Linear Algebra and its Applications}, 746:165--204, 2026.

\bibitem[PSW21]{PlaumannSinnWeis2021}
Daniel Plaumann, Rainer Sinn, and Stephan Weis.
\newblock {Kippenhahn}'s theorem for joint numerical ranges and quantum states.
\newblock {\em SIAM Journal on Applied Algebra and Geometry}, 5(1):86--113,
  2021.

\bibitem[SW{\.Z}18]{SzymanskiWeisZyczkowski2018}
Konrad Szyma{\'n}ski, Stephan Weis, and Karol {\.Z}yczkowski.
\newblock Classification of joint numerical ranges of three {Hermitian}
  matrices of size three.
\newblock {\em Linear Algebra and its Applications}, 545:148--173, 2018.

\bibitem[XZC{\etalchar{+}}20]{XieEtAl2020}
Jie Xie, Aonan Zhang, Ningping Cao, Huichao Xu, Kaimin Zheng, Yiu-Tung Poon,
  Nung-Sing Sze, Ping Xu, Bei Zeng, and Lijian Zhang.
\newblock Observing geometry of quantum states in a three-level system.
\newblock {\em Physical Review Letters}, 125:150401, 2020.

\bibitem[ZH08]{ZachlinHochstenbach2008}
Paul~F. Zachlin and Michiel~E. Hochstenbach.
\newblock On the numerical range of a matrix.
\newblock {\em Linear and Multilinear Algebra}, 56(1--2):185--225, 2008.
\newblock English translation of {Rudolf Kippenhahn}'s 1951 paper.

\end{thebibliography}
\endgroup

\bigskip
{\small
\noindent Carlos Rito\\
Centro de Matem\'atica, Universidade do Minho -- Polo CMAT-UTAD\\
Universidade de Tr\'as-os-Montes e Alto Douro, UTAD\\
Quinta de Prados, 5000-801 Vila Real, Portugal\\
\url{https://www.utad.pt}, \texttt{crito@utad.pt}
\par}

\bigskip

\end{document}